\documentclass[11pt]{amsart}
\usepackage{amssymb}

\usepackage{amscd}

\theoremstyle{definition}
\theoremstyle{remark}
\numberwithin{equation}{section}

\theoremstyle{plain}
\newtheorem{theorem}{Theorem}[section]
\newtheorem{corollary}[theorem]{Corollary}
\newtheorem{lemma}[theorem]{Lemma}
\newtheorem{proposition}[theorem]{Proposition}

\newtheorem{definition}{Definition}[section]
\newtheorem{remark}{Remark}[section]

\begin{document}
\title{Cheeger constant and algebraic entropy of linear groups}
\author{Emmanuel Breuillard}
\email{emmanuel.breuillard@math.univ-lille1.fr}
\author{Tsachik Gelander}
\email{tsachik.gelander@yale.edu}
\date{}
\subjclass{Primary 20Fxx; Secondary 53Cxx}
\keywords{Linear groups, growth, free subgroups.}

\begin{abstract}
We prove a uniform version of the Tits alternative. As a consequence, we
obtain uniform lower bounds for the Cheeger constant of Cayley grahs of
finitely generated non virtually solvable linear groups in arbitrary
characteristic. Also we show that the algebraic entropy of discrete
subgroups of a given Lie group is uniformly bounded away from zero.
\end{abstract}

\maketitle

In this note, we summarize some results whose full proofs will appear in
\cite{BG2}.

\section{Free subgroups in linear groups}

Let $K$ be an arbitrary field and $\Gamma $ a subgroup of $GL_{d}(K)$
generated by a finite subset $\Sigma $. Assume $\Sigma $ is symmetric (i.e. $%
s\in \Sigma \Rightarrow s^{-1}\in \Sigma $), contains the identity $e$, and
let $\mathcal{G}=\mathcal{G}(\Gamma ,\Sigma )$ be the associated Cayley
graph. The set $\Sigma ^{n}$ is the set of all products of at most $n$
elements from $\Sigma ,$ i.e. the ball of radius $n$ centered at the
identity in $\mathcal{G}$. We introduce the following definition:

\begin{definition}
Two elements in a group are said to be \textbf{independent} if they generate
a non-commutative free subgroup. The \textbf{independence diameter} of a
Cayley graph $\mathcal{G}(\Gamma ,\Sigma )$ is the quantity $d_{\Gamma
}(\Sigma )=\inf \{n\in \Bbb{N}$, $\Sigma ^{n}$ contains two independent
elements $\}.$ Similarly, we define the independence diameter of the group $%
\Gamma $ to be $d_{\Gamma }=\sup \{d_{\Gamma }(\Sigma ),$ $\Sigma $ finite
symmetric generating set with $e\in \Sigma \}.$
\end{definition}

The Tits alternative \cite{Tits} asserts that either $\Gamma $ is virtually
solvable (i.e. contains a solvable subgroup of finite index) or $\Gamma $
contains two independent elements, i.e. $d_{\Gamma }(\Sigma )<+\infty $
for every generating set $\Sigma $. The two events are mutually exclusive.
Tits' proof provides no estimate as to how close to the identity in $%
\mathcal{G}$ the independent elements may be found. We obtain:

\begin{theorem}
\label{thm1}(Uniform Tits alternative) Let $\Gamma $ be a finitely generated
subgroup of $GL_{n}(K).$ Assume that $\Gamma $ is not virtually solvable.
Then $d_{\Gamma }<+\infty .$
\end{theorem}

This result improves a theorem of A. Eskin, S. Mozes and H. Oh who proved in
\cite{EMO} the analogous statement when free subgroup is replaced by free
semigroup. Although only linear groups in characteristic zero where
considered in \cite{EMO}, our proof of Theorem \ref{thm1} treats the general
case. This relies on a general version of the Eskin-Mozes-Oh theorem which
we proved together with A. Salehi-Golsefidi and will appear in \cite{BGS}.
As an important corollary and original motivation to their main result, the
authors of \cite{EMO} obtained that finitely generated linear groups in
characteristic zero either are virtually nilpotent or have uniform
exponential growth. Theorem \ref{thm1} has a similar corollary which implies
not just uniform exponential growth but a uniform estimate on the Cheeger
constant of the Cayley graph. Before stating it, we recall the following
definitions:

\begin{definition}
For a group $\Gamma $ generated by a finite symmetric set $\Sigma ,$ we
define its \textbf{uniform $\ell ^{2}$-Kazhdan constant} $\kappa _{\Gamma }(\Sigma )$ to be the largest $\varepsilon \geq 0$
such that
\begin{equation*}
\max_{s\in \Sigma }\left\| s\cdot f-f\right\| _{2}\geq \varepsilon \cdot
\left\| f\right\| _{2}
\end{equation*}
for all $f\in \ell ^{2}(\Gamma ).$ Similarly, we let $\kappa _{\Gamma }=\inf
\{\kappa _{\Gamma }(\Sigma ),$ $\Sigma $ finite symmetric generating set$\}.$
\end{definition}

It is easy to see that $\kappa _{F_{k}}>0$ for the free group $F_{k}$ (e.g. see \cite{Sha}).
Another related quantity is the uniform Cheeger constant defined as follows:

\begin{definition}
For a group $\Gamma $ generated by a finite symmetric set $\Sigma ,$ the
Cheeger constant of the Cayley graph $\mathcal{G}(\Gamma ,\Sigma )$ is
defined by
\begin{equation*}
h_{\Gamma }(\Sigma )=\inf \{\frac{\#\partial A}{\#A},A\subseteq \mathcal{G}%
(\Gamma ,\Sigma )\}
\end{equation*}
where $\partial A=A\backslash \cap _{s\in \Sigma }sA$ is the
inner boundary of a subset $A.$ And the \textbf{uniform Cheeger constant} of
the group $\Gamma $ is defined by $h_{\Gamma }=\inf \{h_{\Gamma }(\Sigma ),$
$\Sigma $ finite symmetric generating set$\}.$
\end{definition}

The above quantities are easily seen to satisfy the following relations:
\begin{equation*}
\sqrt{8\cdot h_{\Gamma }}\geq \sqrt{2}\cdot \kappa _{\Gamma }\geq \frac{%
\kappa _{F_{2}}}{d_{\Gamma }}
\end{equation*}

In \cite{Arz} (see also \cite{Sha}), a finitely generated group is called \textit{uniformly
non-amenable} if $h_{\Gamma }>0.$ We thus have:

\begin{corollary}
Let $\Gamma $ be a finitely generated subgroup of $GL_{n}(K).$ Assume $%
\Gamma $ is not amenable. Then $\kappa _{\Gamma }>0$ and $\Gamma $ is
uniformly non-amenable.
\end{corollary}

One should compare this result to \cite{GZ}, where it is shown that many
arithmetic groups do not have a uniform lower bound for the Kazhdan constant
with respect to an arbitrary unitary representation.


\begin{corollary}
Under the same assumptions on $\Gamma $, there is a constant $\varepsilon
=\varepsilon (\Gamma )>0$ such that if $\Sigma $ is a finite symmetric
generating subset of $\Gamma ,$ then
\begin{equation*}
\#\Sigma ^{n}\geq \#\Sigma \cdot (1+\varepsilon )^{n}
\end{equation*}
for all positive integers $n.$
\end{corollary}

As in Eskin, Mozes and Oh's original proof, the proof of Theorem \ref{thm1}
makes use of the theory of arithmetic groups. However, when proving Theorem
\ref{thm1}, we first obtain along the way a clear and short proof of the
result of Eskin, Mozes and Oh which is very geometric in nature and does not
use arithmetic groups (see Section 3.1 below). Arithmeticity is really
required when one wants to find a free group instead of just a free
semi-group because it is then crucial to obtain elements of $\Gamma $ that
play ping-pong with good ``separation properties''. Several new ingredients
are needed in the proof of Theorem \ref{thm1} such as the
Borel--Harish-Chandra theorem, the argument behind Kazhdan--Margulis theorem
and some facts from the geometry of symmetric spaces and Bruhat-Tits
buildings. We shall outline some of these arguments in Section \ref{proofs}.

The general case reduces to the arithmetic one by the following:

\begin{lemma}
\label{spec}For any finitely generated non-virtually solvable linear group $%
\Gamma $ there is a global field $\Bbb{K}$, a finite set of valuations $S$
of $\Bbb{K}$, a simple $\Bbb{K}$ algebraic group $\Bbb{G}$, and a
homomorphism $f:\Gamma \to \Bbb{G}(\Bbb{K})$ whose image is Zariski
dense and lies in $\Bbb{G}(\mathcal{O}_{\Bbb{K}}(S))$.
\end{lemma}

Note that $f$ is not injective in general. Here $\mathcal{O}_{\Bbb{K}}(S)$ is the ring of $S$-integers in the number
field $\Bbb{K}$, and $\Bbb{G}(\mathcal{O}_{\Bbb{K}}(S))$ is the set of
elements in $\Bbb{G}(\Bbb{K})$ whose matrix elements lie in $\mathcal{O}_{%
\Bbb{K}}(S)$ under some fixed faithful $\Bbb{K}$-representation of $\Bbb{G}$
in $\Bbb{SL}_{n}$.

For arithmetic groups we obtain the following stronger result:

\begin{theorem}
\label{thm2}Let $\Bbb{K}$ be a number field, $S$ a finite set of places of $%
\Bbb{K}$ containing all Archimedean ones, and let $\Bbb{G}$ be a simple $%
\Bbb{K}$-algebraic group. Then there exists a constant $m=m(\Bbb{K},S,\Bbb{G)%
}\geq 1$ with the following property. For any symmetric set $\Sigma $ in $%
\Bbb{G}(\mathcal{O}_{\Bbb{K}}(S))$ with $e\in \Sigma $, which generates a
Zariski dense subgroup $\Gamma $ of $\Bbb{G}$, $d_{\Gamma }(\Sigma )\leq m$.
\end{theorem}

Finally, let us also remark that the main result of \cite{BG}, i.e. the
connected case of the topological Tits alternative, can be deduced easily
from \ref{thm2}.

\section{Algebraic entropy and discrete subgroups}

If $\Gamma $ is a group, let $\mathcal{C}$ be the set of all finite (not
necessarily symmetric) subsets $\Sigma $ containing $e$ and generating $%
\Gamma $.

\begin{definition}
Two elements in a group are said to be \textbf{positively independent} if
they generate a free semigroup. The \textbf{diameter of positive
independence }of set $\Sigma $ containing $e$ is the quantity $d^{pi}(\Sigma
)=\inf \{n\in \Bbb{N}$, $\Sigma ^{n}$ contains two positively independent
elements$\}.$ Similarly, the diameter of positive independence of the group $%
\Gamma $ is defined by $d_{\Gamma }^{pi}=\sup \{d^{pi}(\Sigma ),$ $\Sigma
\in \mathcal{C}\}.$
\end{definition}

The next definition is more standard:

\begin{definition}
Assume $\Gamma $ is finitely generated. For $\Sigma $ in $\mathcal{C}$ we
can define the \textbf{algebraic entropy} of the pair $(\Gamma ,\Sigma )$ to
be the quantity $S_{\Gamma }(\Sigma )=\lim \frac{1}{n}\log (\#\Sigma ^{n}).$
Similarly, the algebraic entropy of $\Gamma $ is defined by $S_{\Gamma
}=\inf_{\Sigma \in \mathcal{C}}S_{\Gamma }(\Sigma ).$
\end{definition}

It is easy to see that $S_{\Gamma }(\Sigma )$ is either positive for all $%
\Sigma $ in $\mathcal{C}$ or $0$ for all $\Sigma $ simultaneously.
Accordingly, the group $\Gamma $ is said to have exponential or
sub-exponential growth. If $S_{\Gamma }>0$, then $\Gamma $ is said to have
\textit{uniform exponential growth}. It is a consequence of Tits' proof of
the Tits alternative and some additional simple argument for solvable groups
that for a linear group $\Gamma \leq GL_{n}(K)$ generated by a finite set $%
\Sigma $, either $d^{pi}(\Sigma )<+\infty $ (and $\Gamma $ has exponential
growth) or $\Gamma $ is virtually nilpotent, hence $d^{pi}(\Sigma )=+\infty $
and $\Gamma $ actually has polynomial growth. The latter quantities are
related by the following inequality:
\begin{equation*}
S_{\Gamma }\geq \frac{\log 2}{d_{\Gamma }^{pi}}
\end{equation*}

A. Eskin, S. Mozes and H. Oh proved in \cite{EMO} that if $\Gamma \leq
GL_{n}(K)$ is finitely generated non-virtually nilpotent then $d_{\Gamma
}^{pi}<\infty ,$ hence $S_{\Gamma }>0$. In general, the constant $d_{\Gamma
}^{pi}$ depends strongly on $\Gamma $ (see the paragraph concluding this
section) however, for discrete subgroups of Lie groups, as well as for
non-relatively compact subgroups over non-Archimedean local fields we have
the following uniform result. Note further that $\Sigma $ is not assumed to
be symmetric in the following statement:

\begin{theorem}
\label{thm3bis} For every integer $d\geq 1$, there is a constant $m=m(d)\geq
1$ with the following property. Let $k$ be a local field and $\Sigma \subset
GL_{d}(k)$ a subset which generates a non-virtually nilpotent group. Assume
further that either:

\begin{itemize}
\item  the group $\langle \Sigma \rangle $ is discrete, or

\item  the field $k$ is non-Archimedean and $\langle \Sigma \rangle $ is not
relatively compact.
\end{itemize}

Then $d^{pi}(\Sigma )\leq m(d).$
\end{theorem}

This implies:

\begin{corollary}
\label{thm4}(Entropy Gap for Discrete Subgroups) For any integer $d\geq 1,$
there is a constant $s=s(d)>0$ such that
\begin{equation*}
S_{\Gamma }>s
\end{equation*}
for all non-virtually nilpotent finitely generated discrete subgroups $%
\Gamma $ of $GL_{d}(\Bbb{R})\times GL_{d}(k_{1})\times ...\times
GL_{d}(k_{n})$ for any $n\geq 0$ and any non-archimedean local fields $%
k_{1},...,k_{n}.$
\end{corollary}

We also prove the following uniform statement for linear groups over general
fields:

\begin{theorem}
For any $n$ there is $m=m(n)\geq 1$, such that if $d\leq n$ and $\Bbb{K}$ is
an algebraic extension of degree $[\Bbb{K}:\Bbb{F}]\leq n$ over a purely
transcendental extension $\Bbb{F}$ of the prime field $\Bbb{K}_{0}$, and $%
\Sigma \subset GL_{d}(\Bbb{K})$ generates a non-virtually nilpotent group,
then $d^{pi}(\Sigma )\leq m$.
\end{theorem}

The restrictions on $[K:F]$ and on the dimension $d$ are really necessary.
In fact, even for $d=2,$ it is possible to find a sequence $\Sigma _{n}$ of
finite symmetric sets in $SL_{2}(\overline{\Bbb{Q}})$ such that each $\Sigma
_{n}$ generates a non virtually nilpotent group, although no pair of
elements in $\Sigma _{n}^{n}$ generates a free semigroup (see \cite{Bre}).
Similarly, R. Grigorchuk and P. de la Harpe have exhibited in \cite{GH} a
sequence of finitely generated non-virtually solvable subgroups $\Gamma _{n}$
in $SL_{k_{n}}(\Bbb{Z})$ such that $\lim \inf S_{\Gamma _{n}}=0$ and $%
k_{n}\rightarrow +\infty .$

\section{Proofs\label{proofs}}

In this section we sketch the proofs of the above results. Let $k$ be a
local field.

\subsection{Proximal elements}

Like in Tits' original proof of his alternative, one basic ingredient in all
the results above is the so-called \textit{ping-pong lemma}. This ensures
that if two projective transformations $a$ and $b$ are in a suitable
geometric configuration when acting on the projective space $\Bbb{P}(k^{d}),$
then $a$ and $b$ generate a free semigroup, or a free group. To describe
this geometric configuration, we need the following definition:

\begin{definition}
An element $g\in PGL_{d}(k)$ is called $\varepsilon $\textit{-contracting},
for some $\varepsilon >0$, if there exists a projective hyperplane $H$,
called a \textit{repelling hyperplane}, and a projective point $v$ called an
\textit{attracting point} such that $d(gp,v)\leq \varepsilon $ whenever $%
p\in \Bbb{P}(k^{d})$ satisfies $d(p,H)\geq \varepsilon $. Moreover $g$ is
called $(r,\varepsilon )$\textit{-proximal}, for $r>2\varepsilon $, if it is
$\varepsilon $-contracting for some $H$ and $v$ with $d(H,v)\geq r.$
Finally, $g$ is called $\varepsilon $\textit{-very contracting} (resp. $%
(\varepsilon ,r)$\textit{-very proximal}) if both $g$ and $g^{-1}$ are $%
\varepsilon $-contracting (resp. $(r,\varepsilon )$-proximal).
\end{definition}

The distance $d([x],[y])$ on $\Bbb{P}(k^{d})$ is the standard distance $%
d([x],[y])=\frac{\left\| x\wedge y\right\| }{\left\| x\right\| \left\|
y\right\| }$ where $\left\| \cdot \right\| $ is a Euclidean norm on $k^{n}$
if $k$ is $\Bbb{R}$ or $\Bbb{C}$ and the supermum norm if $k$ is
non-Archimedean. We refer the reader to \cite{BG} for elementary properties
of such projective transformations. We then have:

\begin{lemma}
\label{pingpong}(The ping-pong lemma) Assume that $a$ and $b$ are $%
(r,\varepsilon )$-very proximal transformations, for some $r>2\varepsilon $,
and the attracting points of $a$ and $a^{-1}$ (resp. $b$ and $b^{-1}$) are
at least $r$ apart from the repelling hyperplanes of $b$ and $b^{-1}$ (resp.
$a$ and $a^{-1}$), then $a$ and $b$ generate a free group.
\end{lemma}

Such elements $a$ and $b$ will be called \textit{ping-pong players}. Observe
that the conditions imposed on $a$ and $b$ imply that the attracting points
of $a$ and $b$ \textit{must} be at least $r-2\varepsilon $ apart. This
\textit{separation property} is a crucial difficulty encountered when trying
to find generators of a free group as opposed to generators of a mere free
semigroup. Indeed, if one only needs a free semigroup, then no condition on
the distance between the two attracting points is necessary as the following
version of the ping-pong lemma for semigroups show:

\begin{lemma}
\label{ping}(The ping lemma) Assume $\varepsilon \leq \frac{1}{3}$ and $%
r>4\varepsilon ^{2}$. Let $a$ be an $(r,\varepsilon ^{3})$-proximal
transformation with attracting fixed point $v$ and repelling hyperplane $H$.
Let $b$ be a projective transformation such that $bv\neq v,$ $d(bv,H)\geq
\varepsilon $ and such that the global Lipschitz constant of $b$ on $\Bbb{P}%
(k^{d})$ satisfies $Lip(b)\leq \frac{1}{\varepsilon }.$ Then $a$ and $ba$
generate a free semigroup.
\end{lemma}

To exhibit $\varepsilon $-contracting elements, it is useful to look at the
Cartan decomposition of $SL_{d}(k)$ since the ratio between the highest and
second to highest component in the diagonal part of the decomposition
determines the contraction properties of the transformation $g$ on $\Bbb{P}%
(k^{d})$ (see \cite{BG}). When $g$ is a diagonal matrix, then this ratio
coincides with the ratio between the highest eigenvalue of $g$ to the second
highest, and the attracting point of $g$ will be the direction of the
highest eigenvector. This situation prevails when $g$ is only assumed to be
\textit{quasi-diagonal}, meaning that the size of its operator norm is
comparable to its highest eigenvalue. This is the ideal situation, because
the attracting points of $g$ being eigendirections, we have control upon
them. In general, elements in a generating set $\Sigma $ need not be
simultaneously quasi-diagonal, however the following crucial proposition
says that up to conjugating $\Sigma $ and looking at a bounded power $\Sigma
^{d^{2}}$ it is possible to bound the norm of elements in $\Sigma $ by the
maximal eigenvalue.

For $y\in SL_{d}(k)$, let $\Lambda _{k}(y)=\max \{|\lambda |_{k}$, $\lambda $
eigenvalue of $y\}.$ For a bounded set $\Omega $ in $SL_{d}(k)$, let $%
\Lambda _{k}(\Omega )=\sup \{\Lambda _{k}(y),y\in \Omega \}$ and
\begin{equation*}
E_{k}(\Omega )=\inf_{h\in SL_{d}(k)}\{\left\| h\Omega h^{-1}\right\| \}
\end{equation*}
where $\left\| \Omega \right\| :=\sup_{y\in \Omega }\left\| y\right\| $ and $%
\|~\|$ is the operator norm. The quantity $E_{k}(\Omega )$ is comparable (up
to bounded powers) with the minimal exponential displacement of the set $%
\Omega $ acting on the symmetric space or building associated to $SL_{d}(k).$
Clearly $\Lambda _{k}(\Omega )\leq E_{k}(\Omega ).$

\begin{proposition}
\label{energy}\label{prop1}There is a constant $c=c(d)>0$ such that for any
compact subset $\Omega $ in $SL_{d}(k)$ with $e\in \Omega $ we have
\begin{equation}
\Lambda _{k}(\Omega ^{d^{2}})\geq c\cdot E_{k}(\Omega )  \label{ineq0}
\end{equation}
Furthermore, if $k$ is not Archimedean, the same holds with $c=1$.
\end{proposition}

Proposition \ref{energy} is a strong version of Proposition 8.5 of
\cite{EMO}, while its proof is significantly simpler and shorter.

\subsection{Generation of free semigroups and a proof of the Eskin-Mozes-Oh
theorem}

The above Proposition is the main step towards producing a proximal element
in $\Sigma ^{d^{2}}$. Its proof is a rather simple contrapositive argument.
Together with some elementary properties of projective transformations as
studied in \cite{BG}, it is essentially enough to prove the result of \cite
{EMO}, namely that $d_{\Gamma }^{pi}<+\infty $ for all non virtually
solvable finitely generated linear groups $\Gamma $, hence also uniform
exponential growth for such $\Gamma $'s. The ping-pong pair is obtained as
follows.

Taking the Zariski closure and after moding out by the solvable radical,
 we can assume that the group $\Gamma $ generated by $%
\Sigma $ is Zariski dense in some semisimple algebraic group lying in $%
SL_{d}.$ By Selberg's lemma we can assume that $\Gamma $ is torsion free.
Using Lemma \ref{Bezout} below, we see that, up to taking a bounded power of
$\Sigma $, we may assume that $\Sigma $ contains a non trivial semisimple
element. Then at least one eigenvalue of this element is not a root of
unity, hence there exists a local field $k$ such that $\Gamma \leq SL_{d}(k)$
and $\Lambda _{k}(\Sigma )>1+\delta $ where $\delta >0$ depends only on $%
\Gamma $. Let $\alpha \in \Sigma ^{d^{2}}$ be such that $\Lambda _{k}(\alpha
)=\Lambda _{k}(\Sigma ^{d^{2}})$. By Proposition \ref{energy}, one can
conjugate $\Sigma $ inside $SL_{d}(k)$ so that $\Lambda _{k}(\alpha )\geq
c\cdot \left\| \Sigma \right\| .$ Up to considering a suitable wedge power
representation $V_{i}=\Lambda ^{i}k^{d}$, we may assume that $\Lambda
_{k}(A)/\lambda _{k}(A)\geq \Lambda _{k}(A)^{1/d^{2}}$ where $\lambda _{k}(A)
$ is the maximum modulus of the second highest eigenvalue and $A=\Lambda
^{i}(\alpha )$. After this operation we have $\left\| \Sigma \right\| \leq $
$\left( \frac{\Lambda _{k}(A)}{c}\right) ^{d}.$ Applying Lemma \ref{conj}
below we can conjugate further $\Sigma $ in $SL(V_{i})$ and get that for all
$n$ large enough (so that $(1+\delta )^{n/2d^{2}}>3$ say), $A^{n}$ is a $%
(1,1/\Lambda _{k}(A)^{n/2d^{2}})$-proximal transformation with attracting
fixed point $v$ and repelling hyperplane $H$ and is such that $\Lambda
_{k}(A)=\left\| A\right\| ,$ while $\left\| \Sigma \right\| \leq $ $\left(
\frac{3\Lambda _{k}(A)}{c}\right) ^{3d\cdot \dim ^{2}V_{i}}.$ Since $\Gamma $
is Zariski dense, not all elements from $\Sigma $ can fix $v$. Applying
Lemma \ref{Bezout} again, we may find an element $B$ in some bounded power
of $\Sigma $ such that none of the powers $B^{j}$ for $j=1,...,\dim V_{i}$
fixes $v.$ But, as can be seen from Cayley Hamilton's theorem for instance,
at least one of the $B^{j}$'s must send $v$ at least $\varepsilon $ away
from $H$ where $\varepsilon $ is at least some fixed bounded power of $%
\left\| B\right\| ,$ hence of $\Lambda _{k}(A)$ because $\left\| B\right\|
\leq \left\| \Sigma \right\| \leq $ $\left( \frac{3\Lambda _{k}(A)}{c}%
\right) ^{3d\cdot \dim ^{2}V_{i}}.$ Hence $A^{n}$ and $B^{j}A^{n}$ are
ping-players (i.e. generate a free semigroup) as soon as $n$ is larger than
a fixed constant depending only on $d,$ $c$ and $\delta $. We can apply the
ping lemma \ref{ping}. {$\square $ }

In order to find non-torsion semisimple ping-pong players in a bounded ball $%
\Sigma ^{m}$ we have just used the following lemma from \cite{EMO}:

\begin{lemma}
\label{Bezout}(\cite{EMO}) Let $\Bbb{G}$ be a Zariski connected algebraic
group. Given a closed algebraic subvariety $X\subset \Bbb{G}$, there is an
integer $k=k(\Bbb{G},X)$ such that for any subset $\Sigma \subset \Bbb{G}$
with $e\in \Sigma $ generating a Zariski dense subgroup of $\Bbb{G}$, the
set $\Sigma ^{k}$ is not contained in $X$.
\end{lemma}

Also we made use of the following simple lemma:

\begin{lemma}
Suppose $A\in SL_{d}(k)$ satisfies $\Lambda _{k}(A)\geq 2$\label{conj}$%
\lambda _{k}(A)$ where $\lambda _{k}(A)$ is the modulus of the second
highest eigenvalue of $A.$ Then the top eigenvalue $\lambda _{1}$ belongs to
$k$, $|\lambda _{1}|=\Lambda _{k}(A)$ and there exists $h\in SL_{d}(%
\overline{k})$ with $\left\| h\right\| \leq 3^{d}\left\| A\right\| ^{d^{2}}$
such that the matrix $A^{\prime }=hAh^{-1}$ is such that $A^{\prime
}e_{1}=\lambda _{1}e_{1}$ and $A^{\prime }H\subset H$ where $H=\left\langle
e_{2},...,e_{d}\right\rangle $ and $\left\| A_{|H}^{\prime }\right\| \leq
\lambda _{k}(A).$
\end{lemma}

To obtain the uniformity in Theorem \ref{thm3bis} we make use of the
classical Margulis lemma which can be stated as follows:

\begin{lemma}
(The Margulis lemma) There is a constant $\varepsilon =\varepsilon (d)>0$
such that for every finite set $\Sigma $ in $SL_{d}(\Bbb{R})$, if $\Sigma $
generates a non virtually nilpotent discrete subgroup, then $E_{\Bbb{R}%
}(\Sigma )\geq 1+\varepsilon .$
\end{lemma}

By another compactness argument we finally obtain:

\begin{proposition}
For every $d\in \Bbb{N}^{*}$ there is a constant $C=C(d)>0$ and an integer $%
N=N(d)>0$ such that for any finite subset $\Sigma $ in $SL_{d}(k)$ with $%
e\in \Sigma $ such that $\Sigma $ generates a non-virtually solvable group,
if $E_{k}(\Sigma )>C$, then $d^{pi}(\Sigma )\leq N.$
\end{proposition}

In particular the constant $s$ in Corollary \ref{thm4} depends on the
Margulis constant $\varepsilon (d)$ and on the constant $c(d)$ from
Proposition \ref{prop1}.

The argument sketched above treats the case of non-virtually solvable linear
groups. For virtually solvable non-virtually nilpotent linear groups we
apply a different argument using one dimensional affine representations
instead of projective representations, which is based on some tools
developed in [\cite{BG1} Section 10] and in \cite{Bre}.

\subsection{Separation properties, arithmeticity and free subgroups}

Here we give some hints on the proof of our main result, Theorem
\ref{thm2}. For the sake of simplicity, let us restrict ourselves
to the case of Zariski-dense subgroups of $SL_{d}(\Bbb{Z})$. In the last paragraph of
this section, we shall give some indications about the general
case.

The main difficulty comes from the fact that in order to generate a free
group rather than just a free semigroup, one should construct a very
proximal element rather than just a proximal one. To do that, one needs both
a good control on the norms of the generators and good separations
properties. For the latter we shall need arithmeticity.

Proposition \ref{energy} supplies us, for each given $\Sigma $, with a
conjugating element $h$ in $SL_{d}(\Bbb{R})$ such that the norm of $h\Sigma
h^{-1}$ is bounded in terms of $\Lambda (\Sigma )$. However, conjugating by $%
h$ we ``loose the arithmeticity''. The first part of the proof of Theorem
\ref{thm2} consists in replacing the conjugating element $h$ by some $\gamma
\in SL_{d}(\Bbb{Z})$:

\begin{proposition}
\label{reduction} There are positive constants $c_{1}=c_{1}(d)$ and $%
r_{1}=r_{1}(d)$ such that for any subset $\Sigma $ in $SL_{d}(\Bbb{Z})$ that
generates a Zariski-dense subgroup of $SL_{d}(\Bbb{R})$ there exists $\gamma\in SL_{d}(%
\Bbb{Z})$ such that
\begin{equation*}
\left\| \gamma \Sigma \gamma ^{-1}\right\| \leq c_{1}\cdot E(\Sigma )^{r_{1}}
\end{equation*}
\end{proposition}

The proof of Proposition \ref{reduction} relies on the following
quantitative variant of Kazhdan--Margulis theorem\footnote{%
Note that Kazhdan-Margulis theorem (which states merely the existence of
unipotents) is trivial in our case because our lattice is arithmetic,
however we need the quantitative estimate.} which was suggested to us by
G.A. Margulis:

\begin{lemma}
\label{margulis} There are positive constants $k_{0},l_{0}$ such that for
any $h\in SL_{d}(\Bbb{R})$ the group $hSL_{d}(\Bbb{Z})h^{-1}\leq G$ contains
a non-trivial unipotent $u$ with
\begin{equation*}
\Vert u-1\Vert \leq l_{0}\cdot \Vert \pi (g)\Vert ^{-k_{0}}
\end{equation*}
where $\pi :G\to SL_{d}(\Bbb{R})/SL_{d}(\Bbb{Z})$ is the canonical
projection and $\Vert \pi (g)\Vert =\min_{\gamma \in SL_{d}(\Bbb{Z})}\Vert
g\gamma \Vert $.
\end{lemma}

The main part of the proof of Theorem \ref{thm2} relies on the construction
of a very proximal element in a bounded power of $\gamma \Sigma \gamma ^{-1}$
acting on the projective space of a corresponding wedge power. This is done
in three steps, in the first we construct a proximal element, in the second
a very contracting one and in the third a very proximal one.

Note that, using Lemma \ref{Bezout}, we can find a semisimple torsion free
element inside a bounded power $\Sigma ^{k_{1}}$ for some constant $k_1$.
The eigenvalues of this element are algebraic integers of bounded degree,
hence $\Lambda (\Sigma ^{k_{1}})\geq 1+\varepsilon_1$ for some constant $%
\epsilon_1>0$. Combining this observation with Proposition \ref{reduction}
we may therefore assume $\left\| \Sigma \right\| \leq \Lambda (\Sigma
)^{r_{2}}$ for some constant $r_{2}=r_{2}(d)>0,$ after changing $\Sigma $
into $\gamma \Sigma ^{N}\gamma ^{-1}$ for some power $N=N(d)$.

We now find $\alpha \in \Sigma $ such that $\Lambda (\alpha )=\Lambda
(\Sigma )$ and make $\Gamma $ act on the (irreducible) wedge power
representation $V_{i}=\Lambda ^{i}\Bbb{C}^{d}$, where $i$ is chosen so that $%
\Lambda (A)/\lambda (A)\geq \Lambda (\Sigma )^{1/d}$ and $\lambda (A)$ is
the maximum modulus of the second highest eigenvalue and $A=\Lambda
^{i}(\alpha ).$ Changing $\Sigma $ into its image under this representation,
we get
\begin{equation*}
\left\| \Sigma \right\| \leq \Lambda (A)^{r_{3}}
\end{equation*}
for some other constant $r_{3}=r_{3}(d)>0.$ This produces the desired
proximal element, hence concludes Step 1.

In order to find a very contracting element in a bounded power of $\Sigma $
we need to find an element $B$ with good separation properties with respect
to $A$, namely a $B$ that sends one eigenvector of $A$ far from some
hyperplane spanned by $d-1$ eigenvectors of $A$.

Part of the problem is to make the word ``far'' more explicit. The condition
that a matrix sends a vector outside a given hyperplane (with no condition
on how far) is an algebraic one and is easily fulfilled thanks to Lemma \ref
{Bezout}. If the matrix has integer coefficients, as it will be the case
thanks to Proposition \ref{reduction} above, then it is possible to estimate
this gap in terms of the norm of the matrix and the arithmetic complexity of
the rationally defined hyperplane. This is the content of Lemma \ref{gap}
below. For a vector $u\in \overline{\Bbb{Q}}^{d}\subset \Bbb{C}^{d}$ we
denote $\left\| u\right\| _{m}=\max \{\left\| \sigma (u)\right\| ,\sigma \in
Gal(\overline{\Bbb{Q}},\Bbb{Q})\}.$

\begin{lemma}
\label{gap}Let $u_{1},...,u_{d}$ be $d$ vectors in $\overline{\Bbb{Q}}^{d}$
whose coordinates are algebraic integers of degree at most $n\geq 1,$ and
let $M:=\max_{i}\left\| u_{i}\right\| _{m}.$ If $H$ is the hyperplane
spanned by $u_{2},...,u_{d}$ and $B\in SL_{d}(\Bbb{Z})$, then either $Bu_{1}\in H,
$ or
\begin{equation*}
d([Bu_{1}],[H])\geq \frac{1}{\left\| B\right\| ^{n^{d}}M^{2dn^{d}}}
\end{equation*}
\end{lemma}

By Lemma \ref{Bezout} we can find $B$ (in a bounded power of $\Sigma $)
which is in ``general position'' with respect to the eigenvectors of $A$. We
then prove:

\begin{proposition}
\label{very-contracting}There is $N=N(d)$ such that the following holds.
Given $q\in \Bbb{N}$, there is a constant $m_{0}>0$ and there is a $B\in
\Sigma ^{N}$ such that the element $A^{m_{0}}BA^{-m_{0}}$ is $\Lambda
(A)^{-q}$-very contracting, with both attracting points lying at a distance
at most $\Lambda (A)^{-q}$ from $v$ -- the eigendirection corresponding to
the maximal eigenvalue of $A$.
\end{proposition}

The proof of Lemma \ref{very-contracting} relies on the following
characterization of contracting elements in terms of their Lipschitz
constants.

\begin{lemma}[See Lemma 3.4, Lemma 3.5 and Proposition 3.3 in \cite{BG1}]
\label{contracting-properties} Let $\epsilon \in (0,\frac{1}{4}],~r\in (0,1]$%
. Let $g\in SL_{n}(\Bbb{R})$ and let $k_{g}a_{g}k_{g}^{\prime }$ be a $KAK$
expression for $g$ where $a_{g}=\text{diag}(a_{1}(g),a_{2}(g),\ldots
,a_{n}(g)),~a_{i}(g)\geq a_{i+1}(g)>0$.

\begin{enumerate}
\item  If $a_{2}(g)/a_{1}(g)\leq \epsilon $ then $g$ is $\epsilon /r^{2}$%
-Lipschitz outside the $r$-neighborhood of the repelling hyperplane $\text{%
span}\{k^{\prime }{}_{g}^{-1}(e_{i})\}_{i=1}^{n}$.

\item  If the restriction of $g$ to some open neighborhood $O\subset \Bbb{P}%
^{n-1}(\Bbb{R})$ is $\epsilon $-Lipschitz, then $a_{2}(g)/a_{1}(g)\leq
\epsilon /2$.

\item  If $a_{2}(g)/a_{2}(g)\leq \epsilon ^{2}$ then $g$ is $\epsilon $%
-contracting, and vice versa, if $g$ is $\epsilon $-contracting, then $%
a_{2}(g)/a_{2}(g)\leq 4\epsilon ^{2}$.
\end{enumerate}
\end{lemma}

The third step is then to obtain a very proximal element by multiplying the
very contracting one $A^{m_{0}}BA^{-m_{0}}$ by some bounded word in the
generators. We want to``separate'' the repelling hyperplanes from the
attracting points. Note that we do not have any information on the position
of the repelling hyperplanes of $A^{m_{0}}BA^{-m_{0}}$, but we do have a
good estimate on the position of its attracting points. We find the right
multiplying element using a simple argument based on the pigeon-hole
principle, and the fact that if $d$ arithmetically defined vectors are
linearly independent in $\Bbb{C}^{d}$ then we can bound from below their
maximal distance to an arbitrary hyperplane in terms of their arithmetic
complexities and norms.

In the last part of the proof of Theorem \ref{thm1}, we conjugate our very
proximal element by a suitable bounded word in the generators and obtain a
second very proximal element which plays ping-pong with the first one. The
argument for finding the appropriate conjugating element is quite similar to
the argument for making the very contracting element a very proximal one.

Let us now say some words about the general case, i.e. when $\Gamma $ is
Zariski-dense in $\Bbb{G}(\mathcal{O}_{\Bbb{K}}(S))$. By the
Borel--Harish-Chandra theorem, $\Bbb{G}(\mathcal{O}_{\Bbb{K}}(S))$ is an
arithmetic lattice in some semisimple Lie group $G\leq \prod_{v\in S}SL_{d}(%
\Bbb{K}_{v})$ over a product of local fields. The absolute value on each $%
\Bbb{K}_{v}$ extends uniquely to any algebraic extension. For $%
g=(g_{v})_{v\in S}$ we define $\Vert g\Vert =\max \Vert g_{v}\Vert _{v}$ and
$\Lambda (g)=\max (\Lambda _{\Bbb{K}_{v}}(g_{v}))$. We obtain the general
version of Proposition \ref{reduction} from Proposition \ref{energy} in two
steps. The fist step consists in replacing the conjugating element $h\in
\prod_{v\in S}SL_{d}(\Bbb{K}_{v})$ by some element $g\in G$. To do that we
exploit theorems of Mostow and Landvogt about totally geodesic imbeddings of
symmetric spaces, and some simple geometric argument using orthogonal
projections on convex subsets in CAT(0) spaces. In a second step, we replace
the element $g$ by some element $\gamma \in \Bbb{G}((\mathcal{O}_{\Bbb{K}%
}(S))$. This step is quite simple in the case where $\Bbb{G}$ is
anisotropic, i.e. when $G/\Bbb{G}(\mathcal{O}_{\Bbb{K}}(S))$ is compact,
because then we can pick $\gamma $ at a bounded distance from $g$. In the
isotropic case we use a generalized version of Lemma \ref{margulis}.

The rest of the proof goes along the same lines sketched above. The guiding
idea is that the distance between arithmetically defined geometric objects
is either $0$ or can be bound from below in terms of their arithmetic
complexity.

\begin{remark}
Let us note that the uniform bound for the independence diameter that our
proof gives for a aubgroup of an arithmetic lattice $\Delta =\Bbb{G}(\mathcal{O}_{\Bbb{K}%
}(S))\leq G$ does strongly depend on $\Delta $ and not just on the ambient
Lie group $G$. In case $\Delta $ is a uniform lattice, it depends on the the
diameter and the injectivity radius of the associated locally symmetric
manifold $K\backslash G/\Gamma $. However, for non-uniform arithmetic
lattices $\Delta \leq G$ we do obtain a uniform constant depending only on $G
$.
\end{remark}

\end{document}